\documentclass[a4paper, 11pt]{article}
\usepackage [a4paper,left=2.5cm,bottom=2.5cm,right=2.5cm,top=2.5cm]{geometry}
\usepackage[english]{babel}
\usepackage{amssymb}
\usepackage{amsmath,amsthm, dsfont}
\usepackage{subcaption}
\usepackage{tabularx,array}
\usepackage{graphicx, url}
\usepackage{enumitem} 

\usepackage{mathtools}  
\usepackage{xfrac} 
\usepackage{bbm}

\theoremstyle{plain}
\theoremstyle{plain}\newtheorem{assumption}{}

\newtheorem{theorem}{Theorem}[section]
\newtheorem{lemma}[theorem]{Lemma}

\newtheorem{proposition}[theorem]{Proposition}

\newtheorem{remark}[theorem]{Remark}
\newtheorem{definition}{Definition}

\usepackage{color}
\usepackage[dvipsnames]{xcolor}

\newcommand{\E}{\mathbb{E}}
\newcommand{\R}{\mathbb{R}}

\renewcommand{\P}{\mathbb{P}}

\title{Local asymptotic normality for discretely observed McKean-Vlasov diffusions}

\author{Akram Heidari \thanks{Department of Mathematics, University of Luxembourg, Maison du Nombre, 6 Avenue de la Fonte, 4364 Esch-sur-Alzette,
Luxembourg, Email: akram.heidari@uni.lu} \and Mark Podolskij \thanks{Department of Mathematics, University of Luxembourg, Maison du Nombre, 6 Avenue de la Fonte, 4364 Esch-sur-Alzette,
Luxembourg, Email: mark.podolskij@uni.lu}}

\date{}

\begin{document}
\maketitle

\begin{abstract} 
We study the local asymptotic normality (LAN) property for the likelihood function associated with discretely observed $d$-dimensional McKean-Vlasov stochastic differential equations over a fixed time interval. The model involves a joint parameter in both the drift and diffusion coefficients, introducing challenges due to its dependence on the process distribution. We derive a stochastic expansion of the log-likelihood ratio using Malliavin calculus techniques and establish the LAN property under appropriate conditions.  The main technical challenge arises from the implicit nature of the transition densities, which we address through integration by parts and Gaussian-type bounds. This work extends existing LAN results for interacting particle systems to the mean-field regime, contributing to statistical inference in non-linear stochastic models. \\


\\
\end{abstract}

\section{Introduction} \label{sec1}
\setcounter{equation}{0}
\renewcommand{\theequation}{\thesection.\arabic{equation}}

The study of McKean-Vlasov stochastic differential equations (SDEs) have gained significant attention in recent years due to their wide-ranging applications in statistical physics, finance, and mean-field games among other fields \cite{BFFT12,CDLL19,CJLW17,DGZZ,FS13,GSS19,GH11,ME99}. These equations are characterized by their dependence on the law of the solution, making them inherently nonlinear. In this work, we consider an i.i.d. array of $d$-dimensional processes, defined  on a complete probability space $(\Omega, \mathcal{F},(\mathcal{F}_t)_{t\geq 0},\P )$, governed by the  McKean-Vlasov SDE:
\begin{equation} \label{model lan} 
\begin{cases}
d X_t^{i,\theta} = b_{\theta_1} \big(X_t^{i,\theta}, \mu_t^{\theta} \big) dt + a_{\theta_2} \big(X_t^{i,\theta} \big) d W_t^i  \qquad i = 1, ... , N, \quad t\in [0, T] \\[1.5 ex]
\text{Law} \big( X_0^{1,{\theta} }, ... , X_0^{N,{\theta}} \big) : = \mu_0 \times ... \times \mu_0
\end{cases}
\end{equation}
where the unknown parameter $ \theta := (\theta_1, \theta_2)$ is an element of the set $\Theta = \Theta_1 \times \Theta_2$, and $\Theta_j \subset \R , j=1,2$ are compact and convex sets with nonempty interior. The $d$-dimensional Brownian motions $(W^i)_{1\leq i \leq N}$  are independent, $\mu_t^{\theta}$ denotes the law of $X_t^{i,\theta}$, and
\[
b: \Theta_1 \times \R^d \times \mathcal{P}_2 \mapsto \R^d, \qquad a: \Theta_2 \times \R^d  \mapsto \R^{d\times d}
\]
are the drift and diffusion coefficient, respectively. Here, $\mathcal{P}_2$ denotes the set of probability measures on $\R^d$ with a finite second moment, endowed with the Wasserstein 2-metric
$$ W_2(\mu, \lambda) := \Big( \inf_{m \in \Gamma (\mu, \lambda)} \int_{\R^d\times \R^d} \|x - y\|^2 m(dx, dy) \Big)^{\frac{1}{2}} $$
where $\Gamma (\mu, \lambda)$ denotes the set of probability measures on the product space $\R^d \times \R^d$ with marginals $\mu$ and $\lambda$.

In this paper, we aim to establish the local asymptotic normality (LAN) property for the likelihood function associated with the discrete observations
\begin{equation}
\left(X_{t_{k}}^{i,\theta}\right)_{k= 1, \dots, n}^{i=1, \ldots, N},
\end{equation}
where $t_{k} := {T k/n}$ and $\Delta_n: = {T/n}$ denotes the discretization step. We consider the asymptotic regime 
$\Delta_n \to 0$, ${N \to \infty}$, and the time horizon $T$  being fixed. We recall that a sequence of statistical models $(\P_{m,\theta}: \theta \in \Theta \in \R^p)$ is said to be locally asymptotically normal if there exist matrices $r_m,\Sigma_{\theta} \in \R^{p\times p}$ such that for any $h\in \R^p$:
\[
\log \frac{\P_{m,\theta + r_m^{-1}h}}{\P_{m,\theta }} = h^{\top} \mathcal{N} -\frac 12 h^{\top} \Sigma_{\theta} h + o_{\P_{m,\theta }}(1), 
\qquad \text{as } m\to \infty,
\]
where $\mathcal{N}\sim \mathcal{N}(0, \Sigma_{\theta})$.
The LAN property is a crucial tool in asymptotic statistical inference, originally introduced by Le Cam. When the LAN property holds and the covariance matrix  is invertible, minimax theorems can be applied to derive lower bounds on the asymptotic variance of estimators (see e.g. \cite{Hajek,lecam,cam&yang}).

Our work contributes to the growing field of statistical inference for McKean-Vlasov processes, providing new insights into the asymptotic properties of parameter estimators in a mean-field setting. Numerous parametric estimation methods for 
McKean-Vlasov diffusions have been studied in the literature under different sampling schemes in \cite{First-paper, B11, C21, CG23,  GL21a, GL21b, K90, L22, SKPP21, WWMX16}. Non-parametric approaches have also gained increasing attention, with notable contributions including \cite{ABPPZ23, BPP22, CGL24, Marc, NPR24}.  

A major challenge in establishing the local asymptotic normality property for McKean–Vlasov SDEs lies in the intractability of their transition densities, which complicates the analysis of the likelihood function’s asymptotic behavior. To address this, we employ tools from Malliavin calculus—specifically, the integration by parts formula—to derive an explicit representation of the logarithmic derivative of the transition density. This methodology, originally introduced by Gobet in the context of classical diffusion models \cite{Gobet 2001, Gobet 2002}, enables a stochastic expansion of the log-likelihood ratio and ultimately yields the desired LAN result. In contrast to Gobet’s framework \cite{Gobet 2002}, which relies on the ergodicity of the underlying process, our approach does not require this assumption, as the asymptotic regime is instead driven by the growing number of particles.

We remark that the drift and diffusion coefficients are estimated at different asymptotic rates. Specifically, the drift parameter is estimated at rate $\sqrt{N}$, while the diffusion parameter is estimated at rate $\sqrt{N/\Delta_n}$. These findings are consistent with the recent results in \cite{First-paper}, which develop a contrast-based estimation method for discretely observed particle systems and establish consistency and asymptotic normality of the resulting estimators. A closely related contribution is \cite{Hof2}, which proves the LAN property for drift estimation in continuously observed 
$d$-dimensional McKean–Vlasov models. It is important to note, however, that in the latter setting the likelihood function admits a closed-form expression via the Girsanov theorem, which simplifies the analysis.

The paper is organized as follows. In Section \ref{secNot} we introduce some notations, formulate the assumptions  and state some technical lemmas we will use in order to show our main results. Section \ref{secMai} is devoted to  the preliminary results essential for proving the LAN property, such as an explicit expression for the  derivative of the log-transition density using the Malliavin calculus, and the main theorem of the paper. Proofs are collected in Section \ref{secProofs}. 
\\

\section{Notation and assumptions}\label{secNot}

In this section, we introduce the notation and the main assumptions associated with the model \eqref{model lan}.

\subsection{Notation}

Throughout the paper, we use a generic positive constant denoted by $C$, or $C_q$ when it depends on an external parameter $q$. This constant is independent of $n$, $N$ and $\theta$, and may vary from line to line. 
The dimensions of the state space and the parameter space, as well as the time horizon $T > 0$, are fixed.

All vectors are understood as column vectors. We denote the Euclidean norm on $\R^d$ by $\| \cdot \|$, and for any vector $x \in \R^d$, we define $x^{\otimes 2} := xx^{\top}$; the $r$-th component of $x$ is written as $x_r$ or $x^r$. The trace of a matrix $A\in \R^{d\times d}$ is denoted by $\text{tr}(A)$. For a function $f:\R^d \times \Theta_j \to \R$, we denote by $\nabla_x f$ (resp. $\partial_{\theta_j} f$) the derivative with respect to $x$ (resp. with respect to $\theta_j$).
 A function $f: \R^d \times \mathcal{P}_l \to \R^d$ is said to have \textit{polynomial growth} if there exist constants $k, l \geq 0$ such that
\[
\| f(x, \mu) \| \leq C \left(1 + \|x\|^k + W_2^l(\mu, \delta_0) \right)
\]
for all $(x, \mu) \in \R^d \times \mathcal{P}_l$, where $\mathcal{P}_l$ denotes the space of probability measures on $\R^d$ with finite $l$-th absolute moment. For $l \in [1, \infty)$, the Wasserstein-$l$ distance between two probability measures $\mu$ and $\lambda$ on $\R^d$ is defined as
\[
W_l(\mu, \lambda) := \left( \inf_{m \in \Gamma(\mu, \lambda)} \int_{\R^d \times \R^d} \|x - y\|^l \, m(dx, dy) \right)^{1/l},
\]
where $\Gamma(\mu, \lambda)$ denotes the set of all couplings of $\mu$ and $\lambda$.

We often suppress the dependence on $\theta$ and write $X_t$ (or  $X_t^i$) to denote the observed data. Moreover, we define conditional expectations 
$$\
E_{t,x}^{\theta}[Z]:= \E^{\theta}[Z| X_t^{\theta}=x] \qquad \text{and} \qquad \E_{t}^{\theta}[Z]:= \E^{\theta}[Z| \mathcal{F}_t],
$$ 
and introduce  the conditional mean and covariance matrix of $X_{t_{k+1}}^{i,\theta}$ given the starting point $X_{t_k}^{i,\theta} = x \in \R^d$ by
\begin{align} \label{mvdef}
m^{\theta}_{t_k,t_{k+1}}(x) &:= \E_{t_k,x}^{\theta} \left[X_{t_{k+1}}^{i,\theta}\right] \in \R^d, \\
V^{\theta}_{t_k,t_{k+1}}(x) &:= \E_{t_k,x}^{\theta} \left[\left(X_{t_{k+1}}^{i,\theta} - m^{\theta}_{t_k,t_{k+1}}(x)\right)^{\otimes 2}\right] \in \R^{d \times d}. \nonumber
\end{align}
For $\mathcal{F}_t$-measurable random variables $Y_t^i$, $i=1,\ldots, N$, we will often use the notation 
\begin{align} \label{Rtnot}
Y_t^i = R_t^i(\varepsilon), \qquad \varepsilon= \varepsilon(n,N),
\end{align}
if $ \varepsilon^{-1} Y_t^i$ is uniformly bounded (in $(\theta,i,t,n,N)$) in $L^q(\Omega)$ for all $q\geq 1$, that is
\[
\E^{\theta}[ |\varepsilon^{-1} Y_t^i|^q]^{1/q} \leq C_q. 
\]
For a sequence of random variables $(Y_n)_{n\geq 1}$, we write $Y_n=O_{\mathbb{P}}(\varepsilon_n)$ (resp. $Y_n=o_{\mathbb{P}}(\varepsilon_n)$) when 
$\varepsilon_n^{-1}Y_n$ is stochastically bounded (resp. $\varepsilon_n^{-1}Y_n$ converges in probability to $0$).


\subsection{Model assumptions}
In order to get asymptotic properties of the likelihood ratio, it is necessary to put some additional conditions on the coefficients. We shall work under the following assumptions on \eqref{model lan}:

\begin{assumption} \label{assum1}

The initial distribution $\mu_0$ is sub-Gaussian, i.e. there exists a constant $\sigma>0$ such that
$$\mu_{0} (A) \leq C \Phi_{\sigma} (A), \qquad \text{for any Borel set } A,$$
where $\Phi_{\sigma}$ denotes the distribution function of $\mathcal{N}_d(0,\sigma^2 \mathbb{I}_d)$.
\end{assumption} 

\begin{assumption} \label{assum2}

The functions $b_{\theta_1}$ and $a_{\theta_2}$ are bounded uniformly in $\theta \in \Theta$. 
For all $\theta$ there exists $C>0$ such that for all $(x,\mu), (y,\lambda) \in \R^d \times {\mathcal{P}_2}$,
$$\|b_{\theta_1}( x, \mu) - b_{\theta_1}( y, \lambda)\| \le C ( \|x - y\| + W_2(\mu, \lambda) )$$
$$\|a_{\theta_2}( x) - a_{\theta_2}( y)\| \le C  \|x - y\|$$

\end{assumption} 

\noindent
We remark that Assumption \ref{assum2} guarantees existence and uniqueness of the McKean-Vlasov SDE 
 \eqref{model lan}. 
As for the regularity of the drift 
$ b: \Theta_1 \times \R^d \times \mathcal{P}_2 \to \R^d $
the notion of linear differentiability, commonly used in the literature on  McKean-Vlasov equations and mean-field games in order to quantify the smoothness of $ \mu \mapsto b_{\theta_1} (x, \mu)$ as a mapping
$ \mathcal{P}_2 \to \R $, will be useful in our setting. We refer in particular to \cite[Section 2]{Marc}
and the references therein.
\begin{definition} \label{deriv of mu}
A mapping $ f: \mathcal{P}_2 \to \R^d $ is said to have a linear functional derivative, if there exists $ \partial_{\mu} f : \R^d \times \mathcal{P}_2 \to \R^d $ such that
\begin{equation*}
f(\mu)-f(\mu^\prime) = \int_0^1 \int_{\R^d} \partial_\mu f(y, \lambda \mu + (1-\lambda) \mu^\prime) (\mu - \mu^\prime) (d y) d\lambda
\end{equation*}
for every $(\mu , \mu^\prime) \in \mathcal{P}_2 $ and $\partial_{\mu} f$ satisfies additional smoothness properties, which will be provided in the following assumption.
\end{definition}

\noindent
In the linear case $f(\mu)=\int_{\R^d} g(x) \mu(dx)$, where $g:\R^d \to\R^d$ is a $\mu$-integrable function, we simply have $\partial_{\mu}f(y,\mu)= g(y)$.

\begin{assumption} \label{assum3}
\textit{Regularity of the derivatives}:
\\
\textit{(I) {For all $(x,\mu)$}, the functions $b_{\theta_1}(x,\mu)$,  
$a_{\theta_2}( x)$ are in $C^2(\Theta_1;\R^d)$,
$C^2(\Theta_2;\R^{d\times d})$ respectively. Furthermore, all their {partial} derivatives up to order three have polynomial growth, 
 uniformly in $\theta$.
 \\
(II) The first and second order derivatives in $\theta$ are locally Lipschitz in $(x,\mu)$ with polynomial weights, i.e.\ {for all $\theta$ there exists $C >0$, $k,l =0,1,\dots$ such that for all $r_1 + r_2 =1, 2$
$(x,\mu), (x^\prime,\mu^\prime) \in \R^d \times {\mathcal{P}_2}$,}
\begin{align*}
&\big\|\partial_{\theta_{1}}^{r_1} \partial_{\theta_{1}}^{r_2} b_{\theta_1} ( x, \mu)- \partial_{\theta_{1}}^{r_1} \partial_{\theta_{1}}^{r_2} b_{\theta_1} (x^\prime,\mu^\prime) \big\| + \big\|\partial_{\theta_{2}}^{r_1} \partial_{\theta_{2}}^{r_2} a_{\theta_2} ( x)-\partial_{\theta_{2}}^{r_1} \partial_{\theta_{2}}^{r_2} a_{\theta_2} ( x^\prime) \big\|
\\
&\qquad \le C (\|x-x^\prime\| + W (\mu,\mu^\prime)) \big( 1 + \|x\|^k + \|x^\prime\|^k + W^l(\mu,\delta_0) + W^l(\mu^\prime,\delta_0) \big).
\end{align*}
\\
(III) The map $\mu \mapsto b_{\theta_1} ( x, \mu)$ admits a functional derivative in the sense of Definition \ref{deriv of mu}.  There exists $k,k',l\geq 1$ such that
\begin{align*}
& \big\|\partial_{\mu} b_{\theta_1} ( x,y, \mu)- \partial_{\mu} b_{\theta_1} ( x^\prime, y^\prime,\mu^\prime) \big\| \le C (\|x-x^\prime\| + \|y-y^\prime\| + W (\mu,\mu^\prime)) \\
& 
\big\|\partial_{\mu} b_{\theta_1} ( x,y, \mu) \big \| \le C  \left(1 + \|x\|^k +\|y\|^{k'}+  W_l^l(\mu, \delta_0) \right).
\end{align*}}
\end{assumption}

\begin{assumption} \label{assum4}

For every point $\theta$ in the interior of $\Theta$, the coefficient functions are twice differentiable in $x$. Furthermore, the following estimates hold: 
 \begin{enumerate} [label=(\alph*)]
 \item $\|\nabla_{x} b_{\theta_1}(x,\mu)\|+\|\nabla_{x} a_{\theta_2}(x)\|\le C $;
 \item $\|g(\cdot, x)\|\le C (1+\|x\|^q)$ for $ g= \nabla^2_{x} b, \nabla_{x} \partial_{\theta_1}b,\nabla^2_{x} a, \nabla_{x} \partial_{\theta_2}a $.
\end{enumerate}
 for some positive constants $C$ and $q$.
\end{assumption} 
\begin{assumption} \label{assum5}
(\textit{Regularity of the diffusion coefficient})
The diffusion matrix $a$ is symmetric, positive definite and satisfies an uniform ellipticity condition. That is, there exist some positive constants $c$ such that 
$$ \forall (\theta_2, x)\in \Theta_2  \times \R^d , \quad \frac{1}{c}  \mathbb{I}_d\le a_{\theta_2}(x)\le c \mathbb{I}_d. $$
\end{assumption} 

\noindent
We note that Assumptions \ref{assum3}-\ref{assum5} are direct analogues of conditions imposed in \cite{Gobet 2001, Gobet 2002}, which in particular imply uniform lower and upper bounds for the transition density established in Proposition \ref{bounds}.

\subsection{Basic elements of Malliavin calculus}

In this section, we present the fundamental concepts of Malliavin calculus that will be used throughout the remainder of the text.

Let \( \mathbb{H} := L^2([0, T], \mathbb{R}^d) \) and consider a function \( h \in \mathbb{H} \). We define the Itô integral of \( h \) with respect to a standard \( \mathbb{R}^d \)-valued Brownian motion \( B = (B_t)_{t \in [0,T]} \) by
\[
B(h) := \int_0^T h(t)^\top dB_t.
\]
Let \( \mathcal{S} \) denote the class of \emph{smooth random variables}, i.e., random variables of the form
\[
F = f(B(h_1), \dots, B(h_k)),
\]
where \( h_i \in \mathbb{H} \) and \( f: \mathbb{R}^k \to \mathbb{R} \) is a \( C_p^\infty \)-function (that is, a function which is infinitely differentiable and all its derivatives grow at most polynomially).

The \textit{Malliavin derivative} of \( F \in \mathcal{S} \) is defined as the \( \mathbb{H} \)-valued random variable
\[
\mathcal{D}F = \sum_{j=1}^k \partial_{x_j} f(B(h_1), \dots, B(h_k)) \cdot h_j.
\]
The operator \( \mathcal{D} \) is closable from \( L^p(\Omega) \) to \( L^p(\Omega; \mathbb{H}) \) for any \( p \geq 1 \). Its closure has domain denoted by \( \mathbb{D}^{1,p} \), which consists of all \( F \in L^p(\Omega) \) such that \( \mathcal{D}F \in L^p(\Omega; \mathbb{H}) \), equipped with the norm
\[
\|F\|_{1,p} := \left( \mathbb{E}[|F|^p] + \mathbb{E}[\|\mathcal{D}F\|_{\mathbb{H}}^p] \right)^{1/p}.
\]
We now introduce the \textit{Skorohod integral} \( \delta \), which is defined as the adjoint of the Malliavin derivative. More precisely, \( \delta \) is a linear operator from a subset of \( L^2([0,T] \times \Omega; \mathbb{H}) \) to \( L^2(\Omega) \), characterized as follows:

\begin{itemize}
    \item[(i)] The domain of \( \delta \), denoted \( \mathrm{Dom}(\delta) \), consists of all processes \( u \in L^2([0,T] \times \Omega; \mathbb{H}) \) such that
    \[
    \forall F \in \mathbb{D}^{1,2}, \quad \left| \mathbb{E}\left[ \langle \mathcal{D}F, u \rangle_{\mathbb{H}} \right] \right| \leq c(u) \|F\|_{L^2(\Omega)}.
    \]
    
    \item[(ii)] For \( u \in \mathrm{Dom}(\delta) \), the Skorohod integral \( \delta(u) \in L^2(\Omega) \) is defined by the \textit{integration by parts formula}
    \begin{equation} \label{eq:skorohod-ibp}
    \forall F \in \mathbb{D}^{1,2}, \quad \mathbb{E}[F \, \delta(u)] = \mathbb{E}[\langle \mathcal{D}F, u \rangle_{\mathbb{H}}].
    \end{equation}
\end{itemize}

\noindent
It is well known that the Skorohod integral and the It\^o integral coincide when both integrals exist. 
The following proposition summarizes key properties of the Skorohod integral.

\begin{proposition} \label{Malliprop}
Let \( p > 1 \). Then the following assertions hold:
\begin{enumerate}[label=(\roman*)]
    \item If \( u \in \mathbb{D}^{1,p}(\mathbb{H}) \), then \( u \in \mathrm{Dom}(\delta) \) and
    \begin{equation} \label{eq:skorohod-bound}
    \|\delta(u)\|_p \leq c_p \left( \|u\|_{L^p(\Omega; \mathbb{H})} + \|\mathcal{D}u\|_{L^p(\Omega; \mathbb{H} \otimes \mathbb{H})} \right).
    \end{equation}

       \item Let \( F \in \mathbb{D}^{1,2} \) and \( u \in \mathrm{Dom}(\delta) \) such that \( \mathbb{E}\left[F^2 \int_0^T \|u_t\|^2 dt \right] < \infty \). Then, whenever the right-hand side is square-integrable, we have the product formula
    \begin{equation} \label{eq:skorohod-product}
    \delta(Fu) = F \delta(u) - \langle \mathcal{D}F, u \rangle_{\mathbb{H}}.
    \end{equation}
\end{enumerate}
\end{proposition}

\section{Main results} \label{secMai}

In this section, we present the main results of the paper. We begin by noting that for any \( s < t \), the conditional law of \( X_t^{\theta,i} \) given \( X_s^{\theta,i} = x \) admits a strictly positive transition density (see Proposition \ref{bounds} and \cite[Proposition 1.2]{Gobet 2002}), denoted by \( p^{\theta}(s, t, x, y) \), which is differentiable with respect to the parameter \( \theta \). 
We define the perturbed parameter as
\begin{align} \label{deftheta+}
(\theta^+_1, \theta^+_2) := \left(\theta^0_1 + \frac{u}{\sqrt{N}}, \; \theta^0_2 + \frac{v}{\sqrt{N/\Delta_n}}\right), 
\qquad \theta^+ := (\theta^+_1, \theta^+_2),
\end{align}
and exploit the Markov property of the processes \( (X_t^{\theta,i})_{t \geq 0} \), \( i = 1, \ldots, N \), to express the log-likelihood ratio between the measures \( \mathbb{P}^{\theta^+} \) and \( \mathbb{P}^{\theta^0} \) as
\begin{align} \label{defzthe}
z(\theta^0,\theta^+) := \log \frac{d \P^{\theta^+}}{d \P^{\theta^0}} (X_{t_k})_{k=1,\ldots, n} = 
\sum_{k=1}^n \sum_{i=1}^N \log \left( \frac{p^{\theta^+}}{p^{\theta^0}} \right) \left(t_k, t_{k+1}, X_{t_k}^i, X_{t_{k+1}}^i \right).
\end{align}
From the expansion \eqref{defzthe}, we obtain the identity
\begin{equation} \label{eq32}
z(\theta^0, \theta^+) = \sum_{k=1}^n \sum_{i=1}^N \big( \zeta_{k}^{i,\theta_1} +  \zeta_{k}^{i,\theta_2}\big),
\end{equation}
where the individual components are given by
\begin{equation} \label{eq32-2}
\begin{split}
\zeta_{k}^{i,\theta_1} &:= \frac{u}{\sqrt{N}} \int_0^1 
\frac{\partial_{\theta_1} p^{\theta_1(l),\theta^+_2}}{p^{\theta_1(l),\theta^+_2}}\left(t_k, t_{k+1}, X_{t_k}^i, X_{t_{k+1}}^i \right) \, dl, \\
\zeta_{k}^{i,\theta_2} &:= \frac{v}{\sqrt{N/\Delta_n}} \int_0^1 
\frac{\partial_{\theta_2} p^{\theta_1^0,\theta_2(l)}}{p^{\theta_1^0,\theta_2(l)}} \left(t_k, t_{k+1}, X_{t_k}^i, X_{t_{k+1}}^i \right) \, dl.
\end{split}
\end{equation}
Here, we use the notation
\begin{align} \label{defil}
\theta_1(l) := \theta^0_1 + \frac{l u}{\sqrt{N}}, \qquad 
\theta_2(l) := \theta^0_2 + \frac{l v}{\sqrt{N/\Delta_n}}.
\end{align}

In the next step, we present an expansion formula for the transition density expressions appearing in 
\eqref{eq32-2}. To this end, we introduce the notation ($t\in[0,\Delta_n]$)
\begin{align*}
X_{t_k+t}^{i,\theta} = X_{t_k}^{i,\theta}+ \int_{0}^t b_{\theta_1} \big(X_{t_k+s}^{i,\theta}, \mu_{t_k+s}^{\theta} \big) ds + \int_0^t a_{\theta_2} \big(X_{t_k+s}^{i,\theta} \big) d B_s^i, 
\end{align*}
where the process $B^i$ represents the Brownian motion $W^i$ shifted in time by $t_k$. Furthermore, for $t\in[0,\Delta_n]$
we introduce derivative processes which satisfy the following stochastic differential equations:
\begin{equation}\label{eq25} 
\begin{split}
\partial_{\theta_1} X_t^{i,\theta} &= \int_0^t \Big( \partial_{\theta_1} b_{\theta_1} \big( X_{t_k+s}^{i,\theta}, \mu_{t_k+s}^{\theta} \big) + \nabla_x b_{\theta_1} \big( X_{t_k+s}^{i,\theta}, \mu_{t_k+s}^{\theta} \big) \, \partial_{\theta_1} X_{s}^{i,\theta} \\
&+ \int_{\R^d} \partial_\mu b_{\theta_1} \big( X_{t_k+s}^{i,\theta},y, \mu_{t_k+s}^{\theta} \big) \, \partial_{\theta_1} \mu_{t_k+s}^{\theta}(dy) \Big) \, ds + \sum_{r=1}^d\int_0^t \nabla_x a^r_{\theta_2}\big( X_{t_k+s}^{i,\theta} \big) \, \partial_{\theta_1} X_s^{i,\theta} \, dB_{r,s}^i, \\
\partial_{\theta_2} X_t^{i,\theta} &= \int_0^t \Big( \nabla_x b_{\theta_1} \big( X_{t_k+s}^{i,\theta}, \mu_{t_k+s}^{\theta} \big) \, \partial_{\theta_2} X_s^{i,\theta} +\int_{\R^d} \partial_\mu b_{\theta_1} \big( X_{t_k+s}^{i,\theta},y, \mu_{t_k+s}^{\theta} \big) \, \partial_{\theta_2} \mu_{t_k+s}^{\theta} (dy) \Big) \, ds \\
& + \sum_{r=1}^d \int_0^t \Big( \partial_{\theta_2} a^r_{\theta_2} \big( X_{t_k+s}^{i,\theta} \big) + \nabla_x a^r_{\theta_2} \big( X_{t_k+s}^{i,\theta} \big) \, \partial_{\theta_2} X_s^{i,\theta} \Big) \, dB_{r,s}^i.
\end{split}
\end{equation}
We remark that the third term in the representation of $\partial_{\theta_1} X_t^{i,\theta}$ is indeed finite
due to Assumption (A3)(III) and the fact that the density of $\partial_{\theta_1}\mu_t^{\theta}$ has a sub-Gaussian density due to Proposition \ref{bounds}. The same argument applies to the corresponding term in 
the representation of $\partial_{\theta_2} X_t^{i,\theta}$.

We now present an explicit formula for the Malliavin derivative \(\mathcal{D}_s X_t^{i,\theta} \in \mathbb{R}^{d \times d}\). To this end, we introduce the \(\mathbb{R}^{d \times d}\)-valued stochastic process \((Y_t^{i,\theta})_{t \in [0,\Delta_n]}\), defined by
\begin{align} \label{defYi}
Y_{t}^{i,\theta} := \mathbb{I}_d + \int_{0}^{t} \nabla_x b_{\theta_1} \big(X_{t_k+s}^{i,\theta}, \mu_{t_k+s}^{\theta} \big) Y_s^{i,\theta} \, ds 
+ \sum_{r=1}^d \int_{0}^{t} \nabla_x a_{\theta_2}^r \big( X_{t_k+s}^{i,\theta} \big) Y_s^{i,\theta} \, dB_s^{r,i},
\end{align}
where \(a_{\theta_2}^r\) denotes the \(r\)th column of the matrix-valued function \(a_{\theta_2}\). It can be shown that \(Y_t^{i,\theta}\) is invertible for all \(t \geq 0\), and the Malliavin derivative satisfies the following identity (see \cite[Eq.~(2.59)]{Nualart}):
\begin{align} \label{formuDX}
\mathcal{D}_s X_{t_k+t}^{i,\theta} = Y_t^{i,\theta} (Y_s^{i,\theta})^{-1} a_{\theta_2}(X_{t_k+s}^{i,\theta}) \,  \, \mathbbm{1}_{[0,t]}(s).
\end{align}
Our first result provides an identity for the terms \( \zeta_{k}^{i,\theta_1} \) and \( \zeta_{k}^{i,\theta_2} \) introduced in \eqref{eq32-2}.

\begin{proposition} \label{transform}
Assume that Assumptions~\ref{assum1}--\ref{assum5} hold. Then, for \( x, y \in \mathbb{R}^d \) and \( \theta \in \Theta \), the transition density \( p^{\theta}(t_k, t_{k+1}, x, y) \) is absolutely continuous with respect to \( \theta \), and we have the representation
\begin{equation} \label{eq23}
\begin{split} 
\frac{\partial_{\theta_1} p^{\theta}}{p^{\theta}}(t_k, t_{k+1}, x, y) &= \frac{1}{\Delta_n} \mathbb{E}_{t_k,x}^{\theta} 
\left[ \sum_{r=1}^d \delta \left( \partial_{\theta_1} X_{r,\Delta_n}^{i,\theta} \, U_r^i \right) \bigg| X_{t_{k+1}}^{i,\theta} = y \right], \\
\frac{\partial_{\theta_2} p^{\theta}}{p^{\theta}}(t_k, t_{k+1}, x, y) &= \frac{1}{\Delta_n} \mathbb{E}_{t_k,x}^{\theta} 
\left[ \sum_{r=1}^d \delta \left( \partial_{\theta_2} X_{r,\Delta_n}^{i,\theta} \, U_r^i \right) \bigg| X_{t_{k+1}}^{i,\theta} = y \right],
\end{split}
\end{equation}
where \( U_r^i \) denotes the \( r \)-th column of the \( \mathbb{R}^{d \times d} \)-valued process
\[
U_s^i := \left( \mathcal{D}_s X_{t_{k+1}}^{i,\theta} \right)^{-1} = a_{\theta_2}^{-1}(X_{t_k+s}^{i,\theta}) \, Y_s^{i,\theta} \, \left( Y_{\Delta_n}^{i,\theta} \right)^{-1}, \qquad s \in [0, \Delta_n].
\]
\end{proposition}

\noindent
The application of Malliavin calculus to derive the transformation in \eqref{eq23} is essential for establishing the LAN property of the statistical model. This formula was originally shown by Gobet (see, e.g., \cite[Proposition 4.1]{Gobet 2001}) in the setting of classical stochastic differential equations, and the same reasoning applies in our framework, so we omit the proof. However, in contrast to classical SDEs, the derivative process 
\( \partial_{\theta_1} X^{i,\theta} \) includes an additional functional derivative term \( \partial_{\mu} b_{\theta_1} \), which captures the dependence of the laws \( \mu_t^{\theta} \) on the parameter \( \theta \).

The next result applies Proposition~\ref{transform} to derive an asymptotic expansion of the log-likelihood function \( z(\theta^0, \theta^+) \).

\begin{proposition} \label{expansi}
Assume that Assumptions~\ref{assum1}--\ref{assum5} hold. Define the quantities
\begin{align} \label{zdef}
z^{\theta}_t(x):= \partial_{\theta_1} b_{\theta_1} \big( x, \mu_{t}^{\theta} \big) + \int_{\R^d} \partial_\mu b_{\theta_1} \big( x,y, \mu_{t}^{\theta} \big) \, \partial_{\theta_1} 
\mu_{t}^{\theta}(dy)
\end{align}
and
\begin{align*}
\widehat{\zeta}_k^{i,\theta_1} &:= \frac{u}{\sqrt{N}} \int_0^1 z^{\theta_1(l), \theta_2^+}_{t_k}(X_{t_k}^i)^{\top}
\left[ a_{\theta_2^+}^{-2}(X_{t_k}^i) \left( X_{t_{k+1}}^i - m_{t_k,t_{k+1}}^{\theta_1(l), \theta_2^+}(X_{t_k}^i) \right) \right] dl, \\[1.5ex]
\widehat{\zeta}_k^{i,\theta_2} &:= \frac{v}{\sqrt{N \Delta_n}} \int_0^1 \text{\rm tr}\left[ 
 \partial_{\theta_2} a_{\theta_2(l)}(X_{t_k}^i) \, a_{\theta_2(l)}^{-1}(X_{t_k}^i)  \right. \\[1.5ex]
&\times  \left.  \left( \left( X_{t_{k+1}}^i - m_{t_k,t_{k+1}}^{\theta_1^0, \theta_2(l)}(X_{t_k}^i) \right)^{\otimes 2} - V_{t_k,t_{k+1}}^{\theta_1^0,\theta_2(l)} (X_{t_k}^i)\right) a_{\theta_2(l)}^{-2}(X_{t_k}^i)\right]
 dl.
\end{align*}
Then, we have for $j=1,2$:
\[
\sum_{k=1}^n \sum_{i=1}^N \left( \zeta_k^{i,\theta_j} - \widehat{\zeta}_k^{i,\theta_j} \right) \xrightarrow{\mathbb{P}^{\theta^0}} 0.
\]
\end{proposition}

\begin{remark}[Interacting particle systems] \rm
A model closely related to \eqref{model lan} is the \textit{interacting particle system} described by
\begin{equation} \label{particl}
\begin{cases}
d X_t^{i,N,\theta} = b_{\theta_1} \big(X_t^{i,N,\theta}, \mu_t^{N} \big) dt + a_{\theta_2} \big(X_t^{i,N,\theta} \big) d W_t^i,  \qquad i = 1, \ldots, N,\quad t \in [0, T], \\[1.5 ex]
\text{Law} \big( X_0^{1,N,\theta}, \ldots, X_0^{N,N,\theta} \big) := \mu_0^{\otimes N},
\end{cases}
\end{equation}
where $\mu_t^N := \frac{1}{N} \sum_{i=1}^N \delta_{X_t^{i,N,\theta}}$ denotes the \emph{empirical measure} of the particle system $(X_t^{1,N,\theta}, \ldots, X_t^{N,N,\theta})$. Parametric estimation techniques for such interacting particle systems have been studied, for instance, in \cite{First-paper,Hof2,SKPP21}. In many cases, statistical inference methods developed for the particle system also extend to its mean-field limit \eqref{model lan}. 

However, establishing the LAN property for the interacting particle system under \emph{discrete-time observations} proves to be significantly more challenging, and no general methods are currently available in the literature. A central difficulty lies in the derivation of lower and upper bounds for the transition densities, which are instrumental in proving Proposition~\ref{expansi}---a key step in the LAN analysis.

In the case of model \eqref{model lan}, where the particles are i.i.d., it suffices to obtain bounds for the $d$-dimensional transition density $p^{\theta}(t,s,x,y)$ (cf. Proposition~\ref{bounds}). In contrast, for interacting particle systems, one cannot reduce the analysis to marginal densities, and must instead handle the full $dN$-dimensional joint transition density $\mathbf{p}^{\theta}(t,s,x,y)$. Unfortunately, such bounds in high dimensions are not yet available in the literature, posing a major obstacle in the LAN analysis for the particle system.

We finally note that these mathematical challenges are absent in the setting of \emph{continuous observations} $(X_t^{i,N,\theta})_{t \in [0,T]}$, where the LAN property for the drift component has been successfully established in \cite{Hof2}.
\end{remark}

\noindent
Using the expansion in Proposition~\ref{expansi}, we now derive the LAN property for the statistical  model~\eqref{model lan}.

\begin{theorem} \label{lan-thrm}
Assume that Assumptions~\ref{assum1}--\ref{assum5} hold. Then,
\begin{equation}
z(\theta^0, \theta^+) \xrightarrow{\mathbb{P}^{\theta^0}-\text{\rm law}} 
\begin{pmatrix}
u\\
v
\end{pmatrix}^{\top}
\mathcal{N}^{\theta^0} - \frac{1}{2}
\begin{pmatrix}
u\\
v
\end{pmatrix}^\top
\Sigma^{\theta^0}
\begin{pmatrix}
u\\
v
\end{pmatrix},
\end{equation}
where the matrix \( \Sigma^{\theta^0}\in \R^{2\times 2} \) is given by
\[
\Sigma^{\theta^0} =
\begin{pmatrix}
\Sigma_b^{\theta^0} & 0 \\
0 & \Sigma_a^{\theta^0}
\end{pmatrix},
\]
with
\begin{align*}
\Sigma_b^{\theta^0} &= \int_0^T \int_{\mathbb{R}^d} z_s^{\theta^0}(x)^{\top} 
a_{\theta_2^0}^{-2}(x) z_s^{\theta^0}(x) \, \mu_s^{\theta^0}(dx) \, ds, \\
\Sigma_a^{\theta^0} &= 2 \int_0^T \int_{\mathbb{R}^d} \text{\rm tr} \left( \partial_{\theta_2} a_{\theta_2^0}(x)  a_{\theta_2^0}^{-1}(x)
 \partial_{\theta_2} a_{\theta_2^0}(x)  a_{\theta_2^0}^{-1}(x) \right) \, \mu_s^{\theta^0}(dx) \, ds.
\end{align*}
\end{theorem}

\noindent
Let us now comment on the implications of Theorem~\ref{lan-thrm}. The article \cite{First-paper} investigates parametric estimation of the drift and diffusion coefficients in one-dimensional particle systems of the form \eqref{particl}, observed at discrete time points. The proposed estimator, based on minimising a contrast function, is shown to be asymptotically normal. Under the additional condition \( N\Delta_n \to 0 \), the asymptotic covariance matrix takes the same form as \( \Sigma^{\theta^0} \). 

However, a key distinction arises when comparing the particle system model \eqref{particl} to the i.i.d.\ observation scheme of the McKean--Vlasov system \eqref{model lan}. In the particle system setting, the linear functional derivative \( \partial_{\mu} b_{\theta_1} \) does not appear in the asymptotic covariance. Consequently, the quantity \( z_t^{\theta}(x) \) is replaced by the simpler expression \( \partial_{\theta_1} b_{\theta_1}(x, \mu_t^{\theta}) \). This structural difference has already been highlighted in \cite{Hof2} in the context of continuously observed particle systems.

For ease of exposition, we have considered the case of univariate parameters. Nevertheless, as shown in \cite{Gobet 2002} for ergodic SDEs, the LAN property extends naturally to the multi-parameter setting. While the notation becomes much more cumbersome, the mathematical complexity does not increase significantly. 
Assume, for instance, that \( \theta_1\in \Theta_1 \subset \mathbb{R}^p \) and \( \theta_2 \in \Theta_2 \subset \mathbb{R}^q \). Then the LAN property holds with asymptotic covariance matrix
\[
\Sigma^{\theta^0} =
\begin{pmatrix}
\Sigma_b^{\theta^0} & 0 \\
0 & \Sigma_a^{\theta^0}
\end{pmatrix}
\in \mathbb{R}^{(p+q) \times (p+q)},
\]
where \( \Sigma_b^{\theta^0} \in \mathbb{R}^{p \times p} \) and \( \Sigma_a^{\theta^0} \in \mathbb{R}^{q \times q} \) are defined by
\begin{align*}
\big(\Sigma_b^{\theta^0}\big)_{kl} &= \int_0^T \int_{\mathbb{R}^d} z_{k,s}^{\theta^0}(x)^{\top} 
a_{\theta_2^0}^{-2}(x) z_{l,s}^{\theta^0}(x) \, \mu_s^{\theta^0}(dx) \, ds, \\
\big(\Sigma_a^{\theta^0}\big)_{kl} &= 2 \int_0^T \int_{\mathbb{R}^d} \mathrm{tr} \left( \partial_{\theta_{k,2}} a_{\theta_2^0}(x)  a_{\theta_2^0}^{-1}(x)
 \partial_{\theta_{l,2}} a_{\theta_2^0}(x)  a_{\theta_2^0}^{-1}(x) \right) \, \mu_s^{\theta^0}(dx) \, ds,
\end{align*}
and
\[
z^{\theta}_{k,t}(x) := \partial_{\theta_{k,1}} b_{\theta_1} \big( x, \mu_{t}^{\theta} \big) + \int_{\mathbb{R}^d} \partial_\mu b_{\theta_1} \big( x, y, \mu_{t}^{\theta} \big) \, \partial_{\theta_{k,1}} \mu_{t}^{\theta}(dy).
\]
We omit the full details of the multi-parameter extension, referring the interested reader to the cited literature for exposition of technicalities.

\section{Proofs} \label{secProofs}

\subsection{Preliminary results} \label{secPrelimi}

In this section, we collect several auxiliary results that will be instrumental in establishing the main theorems of this paper.

To begin with, under Assumptions~\ref{assum3} and~\ref{assum4}, the random variables $X_t^{i,\theta}$, $Y_t^{i,\theta}$, $(Y_t^{i,\theta})^{-1}$, $\partial_{\theta_1} X_t^{i,\theta}$, and $\partial_{\theta_2} X_t^{i,\theta}$ belong to the Malliavin space $\mathbb{D}^{1,p}$ for any $t \in [0,T]$ and $p \ge 1$ (see \cite[Section 2.2]{Nualart}). Moreover, the following uniform estimates hold:
\begin{equation} \label{eq26}
 \mathbb{E}^\theta \left[\sup_{0 \le t \le T} \|Z_t\|^p\right] + \sup_{r \in [0,T]} \mathbb{E}^\theta \left[\sup_{r \le t \le T} \|\mathcal{D}_r Z_t\|^p \right] \leq C,
\end{equation}
for $Z_t \in \{X_t^{i,\theta}, Y_t^{i,\theta}, (Y_t^{i,\theta})^{-1}\}$ and $t \in [0,T]$.

Next, we recall a collection of moment bounds adapted from \cite{First-paper}.

\begin{lemma} \label{l: momentss}
Assume that Assumptions~\ref{assum1}–\ref{assum2} hold. Then, for all $p \ge 1$, $0 \le s < t \le T$ with $t-s \le 1$, $i \in \{1, \dots, N\}$, and $N \in \mathbb{N}$, the following estimates are satisfied:
\begin{enumerate}
    \item[(i)] $\sup_{t \in [0,T]} \mathbb{E}[|X_t^{i,\theta,N}|^p] < C$, and moreover, $\sup_{t \in [0,T]} \mathbb{E}[W_p^q(\mu_t^{\theta}, \delta_0)] < C$ for $p \le q$.
    
    \item[(ii)] $\mathbb{E}[|X_t^{i,N} - X_s^{i,N}|^p] \le C (t - s)^{p/2}$.
    
    \item[(iii)] $\mathbb{E}[W_2^p(\mu_t^{\theta}, \mu_s^{\theta})] \le C (t - s)^{p/2}$.
    
\end{enumerate}
\end{lemma}

\smallskip

\noindent 
In the next step, we present a technique which reduces the McKean-Vlasov system introduced in \eqref{model lan} to standard SDEs with time-varying coefficients. Consider the SDE 
\[
d \widetilde{X}_t^{\theta} = b_{\theta_1} \big(\widetilde{X}_t^{\theta}, \mu_t^{\theta} \big) dt + a_{\theta_2} \big(\widetilde{X}_t^{\theta} \big) d W_t,
\qquad  \widetilde{X}_0\sim \mu_0,
\]
obtained by freezing the laws $(\mu_t)_{t\in [0,T]}$ in the original model \eqref{model lan}. Under Assumption \ref{assum2}, the stochastic process $ \widetilde{X}^{\theta}$ follows an SDE with bounded coefficients.  We also note that the processes $\widetilde{X}^{\theta}$ and $X^{\theta,i}$ have the same law due to uniqueness of the solution of \eqref{model lan}. Consequently,  we can transfer certain results from standard SDE setting to McKean-Vlasov
diffusion. In particular, we obtain the following  upper and lower bounds of Aronson type for the transition density $p^{\theta}(t_k, t_{k+1}, x, y)$ (see \cite[Proposition 1.2]{Gobet 2002}).

\begin{proposition} \label{bounds}
Assume that Assumptions~\ref{assum2}–\ref{assum5} are satisfied. Then there exist constants $c > 1$ and $L > 1$ such that for all $(x,y) \in \mathbb{R}^d \times \mathbb{R}^d$ and $t \in [0,T]$, the transition density $p^{\theta}(t_k, t_{k+1}, x, y)$ satisfies the following Gaussian-type bounds:
\begin{align}
p^{\theta}(t_k, t_{k+1}, x, y) &\le \frac{L}{\Delta_n^{d/2}} \exp\left(-\frac{\|x - y\|^2}{c \Delta_n}\right) \exp(c\Delta_n \|x\|^2), \label{P-bound1} \\
p^{\theta}(t_k, t_{k+1}, x, y) &\ge \frac{1}{L \Delta_n^{d/2}} \exp\left(-c \frac{\|x - y\|^2}{\Delta_n}\right) \exp(-c \Delta_n \|x\|^2). \label{P-bound2}
\end{align}
Moreover, for any $m > 1$, there exist constants $c > 1$, $L > 1$, and $q > 0$ such that
\begin{align}
\E_{t_k,x}^{\bar{\theta}} \left[\left| \frac{\partial_{\theta_1} p^{\theta}}{p^{\theta}}(t_k, t_{k+1}, x, X_{t_{k+1}}^i) \right|^m\right] &\le L \Delta_n^{m/2} \exp(c\Delta_n \|x\|^2)(1 + \|x\|)^q, \label{E-bound1} \\
\E_{t_k,x}^{\bar{\theta}} \left[\left| \frac{\partial_{\theta_2} p^{\theta}}{p^{\theta}}(t_k, t_{k+1}, x, X_{t_{k+1}}^i) \right|^m\right] &\le L \exp(c\Delta_n \|x\|^2)(1 + \|x\|)^q. \label{E-bound2}
\end{align}
\end{proposition}

\noindent
As a consequence of Proposition \ref{bounds} and Assumption  \ref{assum1}, we conclude that all probability measures $(\mu_t)_{t\in [0,T]}$
are sub-Gaussian (cf. \cite[Lemma 3.1]{BPZ25}) and we conclude that 
\begin{align} \label{expMome}
\E^{\theta}\left[\exp(r_0 \|X_t^{\theta,i}\|^2) \right] \leq C \qquad \text{for all } t\in [0,T], 
\end{align}
for some $r_0>0$.

The following lemma will also be used throughout this section. Its proof is a direct consequence of the weak law of large numbers and integrability properties, and is therefore omitted.

\begin{lemma} \label{l: Riemann}
Assume that Assumptions~\ref{assum1}–\ref{assum2} are satisfied. Let $f : \mathbb{R}^d \times \mathcal{P}_l \to \mathbb{R}^p$ be a function such that for some $C>0$, $k, l \in \mathbb{N}_0$, and all $(x,\mu), (y,\lambda) \in \mathbb{R}^d \times \mathcal{P}_l$, the following inequality holds:
\begin{equation}\label{cond:Riemannn}
\|f(x,\mu) \| \le C  \left(1 + \|x\|^k +  W_l^l(\mu, \delta_0) \right).
\end{equation}
Assume further that the map $(x,t) \mapsto f(x, \mu_t)$ is integrable with respect to $\mu_t(dx)\,dt$ on $\mathbb{R}^d \times [0,T]$. Then, as $n,N \to \infty$,
\[
\frac{\Delta_n}{N} \sum_{i=1}^N \sum_{k=1}^n f(X_{t_{k}}^{i,\theta}, \mu_{t_k}^{\theta}) \xrightarrow{\mathbb{P}^{\theta}} \int_0^T \int_{\mathbb{R}^d} f(x, \mu_t^{\theta})\, \mu_t^{\theta}(dx)\, dt.
\]
\end{lemma}

\noindent
The following proposition will be frequently used to establish the negligibility of certain stochastic terms in various settings.

\begin{proposition} \label{prop-neg-drift}
Suppose that Assumptions \ref{assum1}--\ref{assum5} hold. Let $H_{t_{k+1}}^i$, $i=1,\ldots, N$, be a sequence of independent $\mathcal{F}_{t_{k+1}}$-measurable random variables. Recall the definition 
\eqref{Rtnot}.
\begin{itemize}
\item[(i)] Assume that for any $\tau>1$:
$$\E_{t_k,x}^{\theta_1(l), \theta_2^+}[H_{t_{k+1}}^i]=0 \qquad \text{and} \qquad \left(\E_{t_k,x}^{\theta_1(l), \theta_2^+}|H_{t_{k+1}}^i|^\tau\right)^{1/\tau}=R_{t_k}^{i}(\Delta_n^{2}).$$
Then we have 
\begin{equation*}
\frac{1}{\Delta_n\sqrt{N}} \sum_{k=1}^n \sum_{i=1}^N \int_0^1  \E_{t_k,X_{t_{k}}}^{\theta_1(l), \theta_2^+} \left[ H^i|X_{t_{k+1}}^{i,\theta_1(l),\theta_2^+}=X_{t_{k+1}}^i \right] dl \xrightarrow{\P^{\theta^0}} 0.
\end{equation*}
\item[(ii)] Assume that for any $\tau>1$:
$$\E_{t_k,x}^{\theta_1^0, \theta_2(l)}[H_{t_{k+1}}^i]=0 \qquad and \qquad \left(\E_{t_k,x}^{\theta_1^0, \theta_2(l)} |H_{t_{k+1}}^i|^\tau\right)^{1/\tau}=R_{t_k}^{i}(\Delta_n^{3/2}).$$
Then we have 
\begin{equation*}
\frac{1}{\sqrt{N\Delta_n}} \sum_{k=1}^n \sum_{i=1}^N \int_0^1  \E_{t_k,X_{t_{k}}}^{\theta_1^0, \theta_2(l)} \left[ H^i_{t_{k+1}}|X_{t_{k+1}}^{i,\theta_1^0,\theta_2(l)}=X_{t_{k+1}}^i \right] dl \xrightarrow{\P^{\theta^0}} 0 .
\end{equation*}
\end{itemize}
\end{proposition}

\noindent
We note that the main challenge in Proposition~\ref{prop-neg-drift} arises from the need to establish convergence under the original measure $\mathbb{P}^{\theta^0}$, which differs from the perturbed measures $\mathbb{P}^{\theta_1(l), \theta_2^+}$ and $\mathbb{P}^{\theta_1^0, \theta_2(l)}$. The proof presented below closely follows the methodology developed in~\cite{Gobet 2002}.

\begin{proof}
We only show part (i) of Proposition \ref{prop-neg-drift} as the second part is shown similarly (cf. \cite{Gobet 2002}). We define the random variable
\[
\xi_k^{N}:= \frac{1}{\Delta_n\sqrt{N}}  \sum_{i=1}^N \int_0^1  \E_{t_k,X_{t_{k}}}^{\theta_1(l), \theta_2^+} \left[ H^i_{t_{k+1}}|X_{t_{k+1}}^{i,\theta_1(l),\theta_2^+}=X_{t_{k+1}}^i \right] dl.
\]
By martingale methods it suffices to prove the following convergence results:
\[
\sum_{k=1}^n \E^{\theta^0}_{t_k}[\xi_k^{N}] \xrightarrow{\P^{\theta^0}} 0 
\qquad \text{and} \qquad \sum_{k=1}^n \E^{\theta^0}_{t_k}[(\xi_k^{N})^2] \xrightarrow{\P^{\theta^0}} 0. 
\]
We start by handling the conditional expectation. We first observe the identity
\begin{align*}
 \sum_{k=1}^n \E^{\theta^0}_{t_k}[\xi_k^{N}] =  \frac{1}{\Delta_n\sqrt{N}} \sum_{k=1}^n \sum_{i=1}^N \int_0^1  \E_{t_k,X_{t_{k}}^i}^{\theta_1(l), \theta_2^+} \left[ H^i_{t_{k+1}} \frac{p^{\theta^0}}{p^{\theta_1(l),\theta^+_2}}(t_k, t_{k+1}, X_{t_k}^i, X_{t_{k+1}}^i) \right] dl.
\end{align*}
Consequently, we deduce the decomposition
\begin{align*}
 \sum_{k=1}^n \E^{\theta^0}_{t_k}[\xi_k^{N}] =: V_{n,N}^1+ V_{n,N}^2,
\end{align*}
with
\begin{align*}
V_{n,N}^1&:=  \frac{1}{\Delta_n\sqrt{N}} \sum_{k=1}^n \sum_{i=1}^N \int_0^1  \E_{t_k,X_{t_{k}}^i}^{\theta_1(l), \theta_2^+} \left[ H^i_{t_{k+1}} \frac{p^{\theta_1^0,\theta_2^+} - p^{\theta_1(l),\theta_2^+}}{p^{\theta_1(l),\theta^+_2}}(t_k, t_{k+1}, X_{t_k}^i, X_{t_{k+1}}^i) \right] dl, \\[1.5 ex] 
V_{n,N}^2&:=  \frac{1}{\Delta_n\sqrt{N}} \sum_{k=1}^n \sum_{i=1}^N \int_0^1  \E_{t_k,X_{t_{k}}^i}^{\theta_1(l), \theta_2^+} \left[ H^i_{t_{k+1}} \frac{p^{\theta^0} - p^{\theta_1^0,\theta_2^+}}{p^{\theta_1(l),\theta^+_2}}(t_k, t_{k+1}, X_{t_k}^i, X_{t_{k+1}}^i) \right] dl.
\end{align*}
For the first term $V_{n,N}^1$ we obtain the identity
\[
V_{n,N}^1 =  - \frac{u}{\Delta_nN} \sum_{k=1}^n \sum_{i=1}^N \int_0^1 \int_0^1  \E_{t_k,X_{t_{k}}^i}^{\theta_1(l), \theta_2^+} \left[ H^i_{t_{k+1}} \frac{\partial_{\theta_1} p^{\theta_1(r),\theta_2^+}}{ p^{\theta_1(r),\theta_2^+}}\frac{p^{\theta_1(r),\theta_2^+}}{p^{\theta_1(l),\theta^+_2}}(t_k, t_{k+1}, X_{t_k}^i, X_{t_{k+1}}^i) \right] dr dl.
\]
Now, we apply H\"older inequality with conjugates $q_1,q_2,q_3>1$, conditions of Proposition \ref{prop-neg-drift}(i)
and the inequalities of Proposition \ref{bounds}, which imply that
\begin{align*}
&|V_{n,N}^1| \leq \frac{C}{N\sqrt{\Delta_n}} \sum_{k=1}^n \sum_{i=1}^N R_{t_k}^{i}(\Delta_n^{2}) \exp(c\Delta_n\|X_{t_k}^i\|^2/q_2) \\[1.5 ex]
&\times \left( 
\Delta_n^{-d/2} \int_{\R^d} \exp\left(-\frac{q_3\|X_{t_k} - y\|^2}{c\Delta_n}\right) 
\exp\left(\frac{c(q_3-1)\|X_{t_k} - y\|^2}{\Delta_n}\right)  \exp(c(2q_3-1) \Delta_n\|X_{t_k}^i\|^2) dy  
\right)^{1/q_3}.
\end{align*}
We note that the above integral is finite when $(q_3-1)c - q_3/c<0$, which can be achieved by choosing $q_3$ close enough to $1$. We also note 
that, for any given constant $C>0$, $\E[\exp(C\Delta_n \|X_{t_k}^i\|^2)]$ is uniformly bounded due to \eqref{expMome} if $\Delta_n$ is small enough. 
Hence, we conclude that 
\begin{align*}
\E[|V_{n,N}^1|]\leq C\Delta_n^{1/2}. 
\end{align*}
The term $V_{n,N}^2$ is handled in exactly the same fashion and we deduce
\begin{align*}
\E[|V_{n,N}^2|]\leq C\Delta_n^{1/2}. 
\end{align*}
Consequently, we have proved that 
\[
\sum_{k=1}^n \E^{\theta^0}_{t_k}[\xi_k^{N}] \xrightarrow{\P^{\theta^0}} 0 .
\]
Concerning the conditional second moment we readily obtain the formula
\begin{align*}
\sum_{k=1}^n \E^{\theta^0}_{t_k}[(\xi_k^{N})^2] &=  \frac{1}{\Delta_n^2 N}  \sum_{k=1}^n
\sum_{i_1,i_2=1}^N \int_0^1 \int_0^1 \E_{t_k,X_{t_{k}}}^{\theta_1(l_1), \theta_2^+} \left[ H^{i_1}_{t_{k+1}}|X_{t_{k+1}}^{i_1,\theta_1(l_1),\theta_2^+}=X_{t_{k+1}}^{i_1} \right] \\[1.5 ex]
 &\times  \E_{t_k,X_{t_{k}}}^{\theta_1(l_2), \theta_2^+} \left[ H^{i_2}_{t_{k+1}}|X_{t_{k+1}}^{i_2,\theta_1(l_2),\theta_2^+}=X_{t_{k+1}}^{i_2} \right] dl_1dl_2
\end{align*}
Applying Jensen's inequality and proceeding exactly as above, we conclude that
\[
\sum_{k=1}^n \E^{\theta^0}[(\xi_k^{N})^2|\mathcal F_{t_k}]  = O_{\P}(\Delta_n^2).
\]
This completes the proof of Proposition  \ref{prop-neg-drift}.
\end{proof}

\subsection{Proof of Proposition \ref{expansi}}

\subsubsection{The drift term} \label{secdrif}
We start by considering the difference $\zeta_k^{i,\theta_1} - \widehat{\zeta}_k^{i,\theta_1}$. Throughout this subsection
we set for simplicity
\[
\theta=(\theta_1,\theta_2):=(\theta_1(l), \theta_2^+).
\]
According to Proposition \ref{transform} we have the identity
\begin{align*}
\zeta_k^{i,\theta_1} = \frac{u}{\Delta_n\sqrt{N}} \int_0^1  \mathbb{E}_{t_k,X_{t_k}}^{\theta} 
\left[ \sum_{r=1}^d \delta \left( \partial_{\theta_1} X_{r,\Delta_n}^{i,\theta} \, U_r^i \right) \bigg| X_{t_{k+1}}^{i,\theta} = X_{t_{k+1}}^i \right] dl. 
\end{align*}
In the following discussion we will study the decomposition 
\begin{align} \label{firsdeco}
&\delta \left( \partial_{\theta_1} X_{r,\Delta_n}^{i,\theta} \, U_r^i \right) = \Delta_n 
z_{r,\theta_1}^{\theta}\left(X_{t_k}^{i,\theta}\right) \\[1.5 ex]
&\times \left[ a^{-2}_{\theta_2}\left(X_{t_k}^{i,\theta}\right)  \left(
X_{t_{k+1}}^{i,\theta}-m^{\theta}_{t_k,t_{k+1}}\left(X_{t_k}^{i,\theta}\right) \right) \right]_r + H_{t_{k+1}}^i. \nonumber
\end{align}
Here the random variable $H_{t_{k+1}}^i$ is a reminder term, which necessarily satisfies $\E_{t_k,x}^{\theta_1(l), \theta_2^+}[H_{t_{k+1}}^i]=0$ since the other two terms in  \eqref{firsdeco} obviously satisfy this identity. Hence, if we prove that 
\[
 \left(\E_{t_k,x}^{\theta}|H_{t_{k+1}}^i|^\tau\right)^{1/\tau}=R_{t_k}^{i}(\Delta_n^{2}) \qquad \text{for any } \tau>1, 
\]
the proof of Proposition \ref{expansi} is completed for $\zeta_k^{i,\theta_1}$ due to Proposition \ref{prop-neg-drift}(i).

We apply the product formula \eqref{eq:skorohod-product} and conclude the identity
\begin{align*}
\delta \left( \partial_{\theta_1} X_{r,\Delta_n}^{i,\theta} \, U_r^i \right) = 
\partial_{\theta_1} X_{r,\Delta_n}^{i,\theta} \delta\left( U_r^i \right) - \int_{0}^{\Delta_n} 
\big(\mathcal{D}_t \partial_{\theta_1} X_{r,\Delta_n}^{i,\theta}\big)^{\top}  U_{r,t}^i  dt .
\end{align*}
In the next step, we recall the definition $U_s^i = a_{\theta_2}^{-1}(X_{t_k+s}^{i,\theta}) \, Y_s^{i,\theta} \, ( Y_{\Delta_n}^{i,\theta} )^{-1}$, 
$s \in [0, \Delta_n]$, and define its approximation via
\begin{align*}
\widehat{U}_s^i := a_{\theta_2}^{-1}(X_{t_k+s}^{i,\theta}).
\end{align*}
With this notation at hand, we finally obtain the decomposition
\begin{align*}
\delta \left( \partial_{\theta_1} X_{r,\Delta_n}^{i,\theta} \, U_r^i \right) = M^i_{n} + H_{n}^{i,1}+ H_{n}^{i,2}+H_{n}^{i,3}
\end{align*}
with
\begin{align*}
M^i_{n}&:= \Delta_n z_{r,\theta_1}^{\theta}\left(X_{t_k}^{i,\theta}\right) \delta\left( \widehat{U}_r^i \right) \\[1.5 ex]
H_{n}^{i,1}&:= - \int_{0}^{\Delta_n} 
\big(\mathcal{D}_t \partial_{\theta_1} X_{r,\Delta_n}^{i,\theta} \big)^{\top} U_{r,t}^i  dt \\[1.5 ex]
H_{n}^{i,2}&:= \left( \partial_{\theta_1} X_{r,\Delta_n}^{i,\theta} -  \Delta_n z_{r,\theta_1}^{\theta}\left(X_{t_k}^{i,\theta}\right) \right)  \delta\left( \widehat{U}_r^i \right) \\[1.5 ex]
H_{n}^{i,3}&:= \partial_{\theta_1} X_{r,\Delta_n}^{i,\theta} \delta\left( U_r^i -\widehat{U}_r^i \right)
\end{align*}
The term $M^i_n$ represents the main contribution while $H_{n}^{i,j}$, $j=1,2,3$, turn out to be negligible. 

We start with the term $H_{n}^{i,1}$. Due to the formula \eqref{eq25}, we readily obtain the statement
\begin{align} \label{partiest}
\sup_{t\in [0,\Delta_n]} \left| \partial_{\theta_1}X_{r,t}^{i,\theta} \right| = R_{t_{k+1}}^i(\Delta_n).
\end{align}
This implies the same statement for the Malliavin derivative:
\[
 \left\| \mathcal{D}_t \partial_{\theta_1} X_{r,\Delta_n}^{i,\theta}\right\| = R_{t_{k+1}}^i(\Delta_n).
\] 
As a consequence, we deduce that 
\[
H_{n}^{i,1}= R_{t_{k+1}}^i(\Delta_n^2).
\]
Now, we handle the term $H_{n}^{i,2}$. We note that
\begin{align*}
 \delta\left( \widehat{U}_r^i \right) = R_{t_{k+1}}^i (\Delta_n^{1/2}).
\end{align*}
On the other hand, we deduce that 
\begin{align*}
\partial_{\theta_1} X_{r,\Delta_n}^{i,\theta} -  \Delta_n z_{r,\theta_1}^{\theta}\left(X_{t_k}^{i,\theta}\right) = R_{t_{k+1}}^i (\Delta_n^{3/2}),
\end{align*}
which readily implies that $H_{n}^{i,2}= R_{t_{k+1}}^i(\Delta_n^2)$. 

Next, we treat the term $H_{n}^{i,3}$.
Due to \eqref{eq:skorohod-bound}, we conclude that 
\[
\delta\left( U_r^i -\widehat{U}_r^i \right) = R_{t_{k+1}}^i(\Delta_n).
\]
Hence, the application of  \eqref{partiest}, implies that $H_{n}^{i,3}= R_{t_{k+1}}^i(\Delta_n^2)$.

Finally, we treat the main term $M^i_n$. First note that $\delta( \widehat{U}_r^i )$ is a regular It\^o integral with respect to $B^i$. 
For the Brownian motion $B^i$ we obtain the following representation:
\begin{align*}
dB_t^i&= a_{\theta_2}^{-1} \big(X_{t_k+t}^{i,\theta} \big) dX_{t_k+t}^{i,\theta} -  a_{\theta_2}^{-1}\big(X_{t_k+t}^{i,\theta} \big)
b_{\theta_1} \big(X_{t_k+t}^{i,\theta}, \mu_{t_k+t}^{\theta} \big) dt \\[1.5 ex]
&= a_{\theta_2}^{-1} \big(X_{t_k}^{i,\theta} \big) dX_{t_k+t}^{i,\theta} + \left(\mathbb{I}_d - a_{\theta_2}^{-1} \big(X_{t_k}^{i,\theta} \big) 
a_{\theta_2} \big(X_{t_k+t}^{i,\theta} \big) \right) dB^i_t \\[1.5 ex]
&+a_{\theta_2}^{-1} \big(X_{t_k}^{i,\theta} \big) b_{\theta_1} \big(X_{t_k+t}^{i,\theta}, \mu_{t_k+t}^{\theta} \big) dt.
\end{align*}
Using the above representation we conclude that 
\begin{align*}
\delta\left( \widehat{U}_r^i \right) &= \sum_{m=1}^d \int_0^{\Delta_n} \left(a^{-1}_{\theta_2}\right)_{rm} \left(X_{t_k+t}^{i,\theta}\right)dB_{t}^{m,i} \\[1.5 ex]
&= \sum_{m=1}^d \int_0^{\Delta_n} \left(a^{-1}_{\theta_2}\right)_{rm} \left(X_{t_k}^{i,\theta}\right)dB_{t}^{m,i} +R_{t_{k+1}}^i(\Delta_n)
\\[1.5 ex]
&= \sum_{m=1}^d  \left(a^{-2}_{\theta_2}\right)_{rm} \left(X_{t_k}^{i,\theta}\right) \int_0^{\Delta_n} dX_{t_k+t}^{m,i,\theta} +R_{t_{k+1}}^i(\Delta_n)
\\[1.5 ex]
&= \left[ a_{\theta_2}^{-2}(X_{t_k}^{i,\theta}) \left( X_{t_{k+1}}^{i,\theta} - m_{t_k,t_{k+1}}^{\theta}(X_{t_k}^{i,\theta}) \right) \right]_r +R_{t_{k+1}}^i(\Delta_n).
\end{align*}
This completes the proof of Proposition \ref{expansi} for the drift component.

\subsubsection{The diffusion term}
We proceed by considering the difference $\zeta_k^{i,\theta_2} - \widehat{\zeta}_k^{i,\theta_2}$. In this subsection we set for simplicity of notation
\[
\theta=(\theta_1,\theta_2):=(\theta_1^0, \theta_2(l)).
\]
As in Section \ref{secdrif}, according to Proposition \ref{transform} we have the identity
\begin{align*}
\zeta_k^{i,\theta_2} = \frac{v}{\sqrt{\Delta_nN}} \int_0^1  \mathbb{E}_{t_k,X_{t_k}}^{\theta} 
\left[ \sum_{r=1}^d \delta \left( \partial_{\theta_2} X_{r,\Delta_n}^{i,\theta} \, U_r^i \right) \bigg| X_{t_{k+1}}^{i,\theta} = X_{t_{k+1}}^i \right] dl. 
\end{align*}
We consider the following decomposition:
\begin{align} \label{secondeco}
&\delta \left( \partial_{\theta_2} X_{r,\Delta_n}^{i,\theta} \, U_r^i \right) 
 =  \left[\partial_{\theta_2} a_{\theta_2}(X_{t_{k}}^{i,\theta}) a_{\theta_2}^{-1}(X_{t_{k}}^{i,\theta}) 
  \left(
X_{t_{k+1}}^{i,\theta}-m^{\theta}_{t_k,t_{k}}\left(X_{t_k}^{i,\theta}\right) \right)\right]_r
 \nonumber  \\[1.5 ex]
&\times
 \left[  a_{\theta_2}^{-2}(X_{t_{k}}^{i,\theta}) 
  \left(
X_{t_{k+1}}^{i,\theta}-m^{\theta}_{t_k,t_{k+1}}\left(X_{t_k}^{i,\theta}\right) \right)\right]_r
\nonumber
\\[1.5 ex]
&- \left[\partial_{\theta_2} a_{\theta_2}(X_{t_{k}}^{i,\theta}) a_{\theta_2}^{-1}(X_{t_{k}}^{i,\theta}) 
V^{\theta}_{t_k,t_{k+1}} \left( X_{t_{k}}^{i,\theta} \right)  a_{\theta_2}^{-2}(X_{t_{k}}^{i,\theta})  \right]_{r,r} 
+H_{t_{k+1}}^i.
\end{align}
As in the previous section, we only need to show that 
\[
 \left(\E_{t_k,x}^{\theta}|H_{t_{k+1}}^i|^\tau\right)^{1/\tau}=R_{t_k}^{i}(\Delta_n^{3/2}) \qquad \text{for any } \tau>1, 
\]
which completes the proof of Proposition \ref{expansi} is completed for $\zeta_k^{i,\theta_2}$ due to Proposition \ref{prop-neg-drift}(ii).

We apply the product formula \eqref{eq:skorohod-product} and conclude the identity
\begin{align*}
\delta \left( \partial_{\theta_2} X_{r,\Delta_n}^{i,\theta} \, U_r^i \right) = 
 \partial_{\theta_2} X_{r,\Delta_n}^{i,\theta} \delta\left( U_r^i \right) - \int_{0}^{\Delta_n} 
\big(\mathcal{D}_t  \partial_{\theta_2} X_{r,\Delta_n}^{i,\theta}\big)^{\top}   U_{r,t}^i  dt .
\end{align*}
Substituting $U_r^i $ by $\widehat{U}_r^i $ as in the previous proof, we readily deduce the approximation
\begin{align*}
 \partial_{\theta_2} X_{r,\Delta_n}^{i,\theta} = \left[\partial_{\theta_2} a_{\theta_2}(X_{t_{k}}^i) a_{\theta_2}^{-1}(X_{t_{k}}^{i,\theta}) 
  \left(
X_{t_{k+1}}^{i,\theta}-m^{\theta}_{t_k,t_{k}}\left(X_{t_k}^{i,\theta}\right) \right)\right]_r +R_{t_{k+1}}^i (\Delta_n)
\end{align*}
and
\begin{align*}
\delta\left( U_r^i \right) = \left[  a_{\theta_2}^{-2}(X_{t_{k}}^{i,\theta}) 
  \left(
X_{t_{k+1}}^{i,\theta}-m^{\theta}_{t_k,t_{k+1}}\left(X_{t_k}^{i,\theta}\right) \right)\right]_r +R_{t_{k+1}}^i (\Delta_n).
\end{align*}
Furthermore, since $V^{\theta}_{t_k,t_{k+1}} ( X_{t_{k}}^{i,\theta} ) =\Delta_n a^2_{\theta_2}(X_{t_{k}}^{i,\theta} ) + R_{t_{k}}^i (\Delta_n^{3/2})$,
we conclude that 
\begin{align*}
\int_{0}^{\Delta_n} 
\big(\mathcal{D}_t  \partial_{\theta_2} X_{r,\Delta_n}^{i,\theta}  \big)^{\top} U_{r,t}^i  dt 
= \left[\partial_{\theta_2} a_{\theta_2}(X_{t_{k}}^{i,\theta}) a_{\theta_2}^{-1}(X_{t_{k}}^{i,\theta}) 
V^{\theta}_{t_k,t_{k+1}} \left( X_{t_{k}}^{i,\theta} \right)  a_{\theta_2}^{-2}(X_{t_{k}}^{i,\theta})  \right]_{r,r} + R_{t_{k+1}}^i (\Delta_n^{3/2}). 
\end{align*}
This concludes the proof of Proposition  \ref{expansi}.

\subsection{Proof of Theorem \ref{lan-thrm}}
In view of Proposition \ref{expansi} and \cite[Theorem VII-5-2 ]{jacod}, it suffices to check the following conditions:
\begin{align}
& \sum_{k=1}^n \sum_{i=1}^N \E^{\theta^0}_{t_k} \Big[\widehat{\zeta}_k^{i,\theta_1}  \Big] \xrightarrow{\P^{\theta^0}} -\frac{1}{2} u^2 \Sigma_b^{\theta^0} \label{lan-1}\\
& \sum_{k=1}^n \E^{\theta^0}_{t_k} \Big[\Big(\sum_{i=1}^N\widehat{\zeta}_k^{i,\theta_1}\Big)^2 \Big] -  \left(\E^{\theta^0}_{t_k} \Big[\sum_{i=1}^N\widehat{\zeta}_k^{i,\theta_1}  \Big]\right)^2 \xrightarrow{\P^{\theta^0}} u^2\Sigma_b^{\theta^0} \label{lan-2}\\
& \sum_{k=1}^n \E^{\theta^0}_{t_k} \Big[\Big|  \sum_{i=1}^N \widehat{\zeta}_k^{i,\theta_1}\Big|^{4}  \Big] \xrightarrow{\P^{\theta^0}} 0
 \label{lan-3} \\
& \sum_{k=1}^n  \sum_{i=1}^N \E^{\theta^0}_{t_k} \Big[\widehat{\zeta}_k^{i,\theta_2} \Big] \xrightarrow{\P^{\theta^0}} -\frac{1}{2} v^2 \Sigma_a^{\theta^0} \label{lan-4b} \\
& \sum_{k=1}^n \E^{\theta^0}_{t_k} \Big[\Big( \sum_{i=1}^N \widehat{\zeta}_k^{i,\theta_2}\Big)^2  \Big] -  \left(\E^{\theta^0}_{t_k} \Big[ \sum_{i=1}^N\widehat{\zeta}_k^{i,\theta_2}  \Big]\right)^2 \xrightarrow{\P^{\theta^0}}  v^2 \Sigma_a^{\theta^0} \label{lan-5b} \\
& \sum_{k=1}^n \E^{\theta^0}_{t_k} \Big[\Big|  \sum_{i=1}^N \widehat{\zeta}_k^{i,\theta_2}\Big|^{4} \Big] \xrightarrow{\P^{\theta^0}} 0  \label{lan-6b} \\
& \sum_{k=1}^n \E^{\theta^0}_{t_k} \Big[ \sum_{i=1}^N\widehat{\zeta}_k^{i,\theta_1}  \sum_{i=1}^N \widehat{\zeta}_k^{i,\theta_2} \Big] - \E^{\theta^0}_{t_k} \Big[  \sum_{i=1}^N \widehat{\zeta}_k^{i,\theta_1} \Big] \E^{\theta^0}_{t_k} \Big[ \sum_{i=1}^N \widehat{\zeta}_k^{i,\theta_2} \Big] \xrightarrow{\P^{\theta^0}} 0 \label{lan-7b}
\end{align}
\textit{Proof of \eqref{lan-1}:} We first observe that $ m^{\theta}_{t_k,t_{k+1}}(x)= x+ \int_{t_k}^{t_{k+1}} \E^{\theta}_{t_k,x} [b_{\theta_1}(X_{s}^{i, \theta} , \mu_{s}^{\theta})] d s $.
Applying the Taylor expansion and recalling \eqref{zdef}, we deduce that
\begin{align} \label{meandiff}
m^{\theta^0}_{t_k,t_{k+1}}(X_{t_k}^i) - m^{\theta_1(l), \theta_2^+}_{t_k,t_{k+1}}(X_{t_k}^i) &= 
-\frac{lu\Delta_n}{\sqrt{N}} z^{\theta^0}_{t_k}(X_{t_k}^i)  + 
R_{t_{k}}^i(\varepsilon_{n,N}\Delta_n/\sqrt{N})
\end{align}
with $\varepsilon_{n,N}\to 0$.
We introduce the notation 
\begin{align*}
\overline{\zeta}_k^{i,\theta_1} & := \frac{u}{\sqrt{N}} \int_0^1 z^{\theta_1(l),\theta_2^+}_{t_k}(X_{t_k}^i)^{\top}
\left[ a_{\theta_2^+}^{-2}(X_{t_k}^i) \left( X_{t_{k+1}}^i - m_{t_k,t_{k+1}}^{\theta_0}(X_{t_k}^i) \right) \right] dl,
\end{align*}
which is a direct analogue of $\widehat{\zeta}_k^{i,\theta_1}$ with $m^{\theta_1(l), \theta_2^+}$ replaced by $m^{\theta_0}$.  By definition it holds that 
$ \E^{\theta^0}_{t_k}[\overline{\zeta}_k^{i,\theta_1} ]=0$. Hence, we conclude the statement
\begin{align*}
\sum_{k=1}^n \sum_{i=1}^N \E^{\theta^0}_{t_k} \Big[\widehat{\zeta}_k^{i,\theta_1}  \Big] &=
 \frac{u}{\sqrt{N}} \sum_{k=1}^n \sum_{i=1}^N \int_0^1 z^{\theta_1(l),\theta_2^+}_{t_k}(X_{t_k}^i)^{\top} \\[1.5ex]
&\times
\left[ a_{\theta_2^+}^{-2}(X_{t_k}^i) \left( m_{t_k,t_{k+1}}^{\theta^0}(X_{t_k}^i)  - m_{t_k,t_{k+1}}^{\theta_1(l), \theta_2^+}(X_{t_k}^i) \right) \right] dl
\xrightarrow{\P^{\theta^0}} -\frac{1}{2} u^2 \Sigma_b^{\theta^0},
\end{align*} 
which follows directly from \eqref{meandiff} and Lemma \ref{l: Riemann}. \\ 

\noindent
\textit{Proof of \eqref{lan-2}:} We observe that
\begin{align*}
\E^{\theta^0}_{t_k} \Big[\sum_{i=1}^N\widehat{\zeta}_k^{i,\theta_1}  \Big]= 
\E^{\theta^0}_{t_k} \Big[\sum_{i=1}^N(\widehat{\zeta}_k^{i,\theta_1}- \overline{\zeta}_k^{i,\theta_1})  \Big] = R_{t_{k}}^i(\Delta_n)
\end{align*}
due to  \eqref{meandiff}. Hence, the second term in \eqref{lan-2} is asymptotically negligible. As for the first term, we obtain the decomposition
\begin{align*}
\sum_{k=1}^n \E^{\theta^0}_{t_k} \Big[\Big(\sum_{i=1}^N\widehat{\zeta}_k^{i,\theta_1}\Big)^2  \Big] 
= \sum_{k=1}^n \E^{\theta^0}_{t_k} \Big[\sum_{i_1\not= i_2}\widehat{\zeta}_k^{i_1,\theta_1} \widehat{\zeta}_k^{i_2,\theta_1} 
+ \sum_{i=1}^N (\widehat{\zeta}_k^{i,\theta_1})^2   \Big]. 
\end{align*}
Again due to \eqref{meandiff}, we deduce the estimate
\begin{align*}
\E^{\theta^0}_{t_k} \Big[\sum_{i_1\not= i_2}\widehat{\zeta}_k^{i_1,\theta_1} \widehat{\zeta}_k^{i_2,\theta_1}    \Big]=
\E^{\theta^0}_{t_k} \Big[\sum_{i_1\not= i_2}(\widehat{\zeta}_k^{i_1,\theta_1} -\overline{\zeta}_k^{i_1,\theta_1})( \widehat{\zeta}_k^{i_2,\theta_1}-\overline{\zeta}_k^{i_2,\theta_1}) 
  \Big]= R_{t_{k}}^i(\Delta_n^2).
\end{align*}
In other words,
\[
 \sum_{k=1}^n \E^{\theta^0}_{t_k} \Big[\sum_{i_1\not= i_2}\widehat{\zeta}_k^{i_1,\theta_1} \widehat{\zeta}_k^{i_2,\theta_1} 
  \Big] \xrightarrow{\P^{\theta^0}} 0.
\]
Finally, we have that 
\begin{align*}
(\widehat{\zeta}_k^{i,\theta_1})^2 =  \frac{u^2}{N}  \left[ z^{\theta^0}_{t_k}(X_{t_k}^i)^{\top}
a_{\theta_2^0}^{-2}(X_{t_k}^i) \left( X_{t_{k+1}}^i - m_{t_k,t_{k+1}}^{\theta^0}(X_{t_k}^i) \right) \right]^{2}  +R_{t_{k+1}}^{i}\left(\varepsilon_{n,N} \Delta_n/N\right)
\end{align*}
with $\varepsilon_{n,N} \to 0$. Since $V^{\theta}_{t_k,t_{k+1}} ( X_{t_{k}}^{i,\theta} ) =\Delta_n a^2_{\theta_2}(X_{t_{k}}^{i,\theta} ) + R_{t_{k}}^i (\Delta_n^{3/2})$, we conclude that 
\begin{align*}
\sum_{k=1}^n \sum_{i=1}^N \E^{\theta^0}_{t_k} \Big[ (\widehat{\zeta}_k^{i,\theta_1})^2   \Big]
\xrightarrow{\P^{\theta^0}} u^2\Sigma_b^{\theta^0}.  
\end{align*}

\noindent
\textit{Proof of \eqref{lan-3}:} Observe that 
\begin{align*}
\sum_{k=1}^n \E^{\theta^0}_{t_k} \Big[\Big|  \sum_{i=1}^N \widehat{\zeta}_k^{i,\theta_1}\Big|^{4}  \Big] 
&\leq C\sum_{k=1}^n \E^{\theta^0}_{t_k} \Big[\Big|  \sum_{i=1}^N \widehat{\zeta}_k^{i,\theta_1} -
\E^{\theta^0}_{t_k} \big[\widehat{\zeta}_k^{i,\theta_1}   \big] \Big|^{4}  \Big] \\[1.5 ex]
&+C\sum_{k=1}^n \Big|\sum_{i=1}^N  
\E^{\theta^0}_{t_k} \big[\widehat{\zeta}_k^{i,\theta_1}   \big] \Big|^{4} . 
\end{align*}
As in the proof of \eqref{lan-2}, we immediately deduce the convergence
\[
\sum_{k=1}^n \Big|\sum_{i=1}^N  
\E^{\theta^0}_{t_k} \big[\widehat{\zeta}_k^{i,\theta_1}   \big] \Big|^{4} \xrightarrow{\P^{\theta^0}} 0.
\]
Setting $y_k^i:= \widehat{\zeta}_k^{i,\theta_1} -
\E^{\theta^0}_{t_k} [\widehat{\zeta}_k^{i,\theta_1}] $, we obtain the decomposition 
\begin{align*}
\sum_{k=1}^n \E^{\theta^0}_{t_k} \Big[\Big|  \sum_{i=1}^N y_k^i \Big|^{4} \Big] 
= \sum_{k=1}^n \E^{\theta^0}_{t_k} \left[ \sum_{(i_1,i_2,i_3,i_4)\in A_N \cup B_N \cup C_N } y_k^{i_1}\cdots y_k^{i_4}  \right], 
\end{align*}
where 
\begin{align*}
A_N&:=\{\text{all indices are distinct}\},\\
B_N&:=\{\text{the indices contain two pairs } (m_1,m_1) \text{ and } (m_2,m_2) \text{ with } m_1\not= m_2\},\\
C_N&:=\{i_1=i_2=i_3=i_4\}.
\end{align*}
Due to independence and \eqref{meandiff} we have
\[
\sum_{k=1}^n \E^{\theta^0}_{t_k} \left[ \sum_{(i_1,i_2,i_3,i_4)\in A_N } y_k^{i_1}\cdots y_k^{i_4}   \right] =O_{\mathbb{P}^{\theta^0}}(\Delta_n^3).
\] 
We also obtain
\[
\sum_{k=1}^n \E^{\theta^0}_{t_k} \left[ \sum_{(i_1,i_2,i_3,i_4)\in B_N } y_k^{i_1}\cdots y_k^{i_4}  \right] =O_{\mathbb{P}^{\theta^0}}(\Delta_n)
\]
and
\[
\sum_{k=1}^n \E^{\theta^0}_{t_k} \left[ \sum_{(i_1,i_2,i_3,i_4)\in C_N } y_k^{i_1}\cdots y_k^{i_4}  \right] =O_{\mathbb{P}^{\theta^0}}(\Delta_n/N).
\]
This completes the proof of  \eqref{lan-3}. \\

\noindent
\textit{Proof of \eqref{lan-4b}:} Similarly to the proof of  \eqref{lan-1}, we obtain the identity:
\begin{align*}
&\E^{\theta^0}_{t_k} \Big[\widehat{\zeta}_k^{i,\theta_2} \Big] = \frac{v}{\sqrt{N \Delta_n}} \int_0^1 \text{\rm tr}\left[ 
 \partial_{\theta_2} a_{\theta_2(l)}(X_{t_k}^i) \, a_{\theta_2(l)}^{-1}(X_{t_k}^i)  \right. \\[1.5ex]
&\times  \left.  \left( \left( m_{t_k,t_{k+1}}^{\theta^0}(X_{t_k}^i)  - m_{t_k,t_{k+1}}^{\theta_1^0, \theta_2(l)}(X_{t_k}^i) \right)^{\otimes 2} +
\left(V_{t_k,t_{k+1}}^{\theta^0} (X_{t_k}^i) - V_{t_k,t_{k+1}}^{\theta_1^0,\theta_2(l)} (X_{t_k}^i)\right)\right) a_{\theta_2(l)}^{-2}(X_{t_k}^i) \right]
 dl.
\end{align*}
Due to  \eqref{meandiff} the term containing the difference of conditional expectations is negligible, and we only need to study the difference of conditional variances. Via It\^o formula we deduce the decomposition 
\begin{align*}
&\left(V_{t_k,t_{k+1}}^{\theta}(x)\right)_{r_1,r_2} = x_{r_1}x_{r_2} + \int_{t_k}^{t_{k+1}} \E^{\theta}_{t_k,x} \left[ 
a_{\theta_2}^{2}(X_s^{i,\theta})_{r_1,r_2} + b_{r_1,\theta_1}(X_s^{i,\theta}, \mu_s^{\theta}) X_{r_2,s}^{i,\theta}
+ b_{r_2,\theta_1}(X_s^{i,\theta}, \mu_s^{\theta}) X_{r_1,s}^{i,\theta} \right]ds \\[1.5 ex]
&\qquad \qquad- m_{t_k,t_{k+1}}^{\theta}(x)_{r_1} m_{t_k,t_{k+1}}^{\theta}(x)_{r_2} \\[1.5 ex]
&= \int_{t_k}^{t_{k+1}} \E^{\theta}_{t_k,x} \left[ 
a_{\theta_2}^{2}(X_s^{i,\theta})_{r_1,r_2} + b_{r_1,\theta_1}(X_s^{i,\theta}, \mu_s^{\theta}) (X_{r_2,s}^{i,\theta}-x_{r_2})
+ b_{r_2,\theta_1}(X_s^{i,\theta}, \mu_s^{\theta}) (X_{r_1,s}^{i,\theta}-x_{r_1}) \right]ds \\[1.5 ex]
&+ \left(\int_{t_k}^{t_{k+1}} \E^{\theta}_{t_k,x} \left[b_{r_1,\theta_1}(X_s^{i,\theta}, \mu_s^{\theta})  \right]ds\right)
\left(\int_{t_k}^{t_{k+1}} \E^{\theta}_{t_k,x} \left[b_{r_2,\theta_1}(X_s^{i,\theta}, \mu_s^{\theta})  \right]ds\right).
\end{align*}
By mean value theorem we deduce 
\[
V_{t_k,t_{k+1}}^{\theta^0} (x) - V_{t_k,t_{k+1}}^{\theta_1^0,\theta_2(l)} (x) = (\theta_2^0 - \theta_2(l)) \partial_{\theta_2} 
V_{t_k,t_{k+1}}^{\theta_1^0,\widetilde{\theta}_2} (x)
\]
for some $\widetilde{\theta}_2$ satisfying $|\theta_2^0 - \widetilde{\theta}_2|\leq |\theta_2^0 - \theta_2(l)|$.
Consequently, we conclude that 
\begin{align} \label{Vconve}
V_{t_k,t_{k+1}}^{\theta^0} (X_{t_k}^i) - V_{t_k,t_{k+1}}^{\theta_1^0,\theta_2(l)} (X_{t_k}^i) = -\frac{2lv\Delta_n^{3/2}}{\sqrt{N}} 
\partial_{\theta_2} a_{\theta_2^0}(X_{t_k}^i) a_{\theta_2^0}(X_{t_k}^i) +R_{t_{k}}^{i}\left(\varepsilon_{n,N} \Delta_n^{3/2}/\sqrt{N}\right)
\end{align}
with $\varepsilon_{n,N} \to 0$. Hence, \eqref{Vconve} implies the convergence in \eqref{lan-4b}. \\

\noindent
\textit{Proof of \eqref{lan-5b}:} In the previous proof we have shown that 
\begin{align*}
\E^{\theta^0}_{t_k} \Big[\widehat{\zeta}_k^{i,\theta_2} \Big] = R_{t_{k}}^i(\Delta_n/N).
\end{align*}
Hence, we immediately conclude that 
\[
\sum_{k=1}^n  \left(\E^{\theta^0}_{t_k} \Big[ \sum_{i=1}^N\widehat{\zeta}_k^{i,\theta_2}  \Big]\right)^2 \xrightarrow{\P^{\theta^0}}  0.
\]
Furthermore, following the same arguments as in the proof of \eqref{lan-2}, we deduce the estimate
\begin{align*}
\sum_{k=1}^n \E^{\theta^0}_{t_k} \Big[\Big(\sum_{i=1}^N\widehat{\zeta}_k^{i,\theta_2}\Big)^2  \Big] 
=  \sum_{i=1}^N \E^{\theta^0}_{t_k} \Big[(\widehat{\zeta}_k^{i,\theta_2})^2   \Big] + o_{\mathbb{P}^{\theta^0}}(1).  
\end{align*}
Last but not least, a direct computation shows that 
\begin{align*}
 \sum_{i=1}^N \E^{\theta^0}_{t_k} \Big[(\widehat{\zeta}_k^{i,\theta_2})^2   \Big] \xrightarrow{\P^{\theta^0}} v^2 \Sigma_a^{\theta^0}.
\end{align*}
This completes the proof of  \eqref{lan-5b}.\\

\noindent
\textit{Proof of \eqref{lan-6b}:} This statement is shown in exactly the same way as  \eqref{lan-3}. \\

\noindent
\textit{Proof of \eqref{lan-7b}:} Applying the estimates from the proof of \eqref{lan-1} and \eqref{lan-4b}, we immediately obtain the convergence
\[
\sum_{k=1}^n  \E^{\theta^0}_{t_k} \Big[  \sum_{i=1}^N \widehat{\zeta}_k^{i,\theta_1} \Big] \E^{\theta^0}_{t_k} \Big[ \sum_{i=1}^N \widehat{\zeta}_k^{i,\theta_2} \Big] \xrightarrow{\P^{\theta^0}} 0. 
\]
Also, using the same arguments as in the proof of \eqref{lan-2}, we deduce that 
\[
\sum_{k=1}^n \E^{\theta^0}_{t_k} \Big[ \sum_{i=1}^N\widehat{\zeta}_k^{i,\theta_1}  \sum_{i=1}^N \widehat{\zeta}_k^{i,\theta_2} \Big]
= \sum_{k=1}^n \E^{\theta^0}_{t_k} \Big[ \sum_{i=1}^N\widehat{\zeta}_k^{i,\theta_1}   \widehat{\zeta}_k^{i,\theta_2} \Big] + o_{\mathbb{P}^{\theta^0}}(1).
\]
Finally, using the representations from Proposition \ref{expansi}, we obtain for any $1\leq r_1,r_2,r_3\leq d$:
\begin{align*}
\E^{\theta^0}_{t_k}\left[\left( X_{t_{k+1}}^i - m_{t_k,t_{k+1}}^{\theta_1(l), \theta_2^+}(X_{t_k}^i) \right)_{r_1}
V_{t_k,t_{k+1}}^{\theta_1^0,\theta_2(l)} (X_{t_k}^i)_{r_2,r_3}\right] = R_{t_k}^i(\Delta_n^2/\sqrt{N}), \\[1.5 ex]
\E^{\theta^0}_{t_k}\left[\left( X_{t_{k+1}}^i - m_{t_k,t_{k+1}}^{\theta_1(l), \theta_2^+}(X_{t_k}^i) \right)_{r_1}
\left( X_{t_{k+1}}^i - m_{t_k,t_{k+1}}^{\theta_1^0, \theta_2(l)}(X_{t_k}^i) \right)^{\otimes 2}_{r_2,r_3}\right] =
R_{t_k}^i(\Delta_n^2)
\end{align*}
This implies the convergence
\[
\sum_{k=1}^n \E^{\theta^0}_{t_k} \Big[ \sum_{i=1}^N\widehat{\zeta}_k^{i,\theta_1}   \widehat{\zeta}_k^{i,\theta_2} \Big] 
\xrightarrow{\P^{\theta^0}} 0,
\]
and the proof of Theorem \ref{lan-thrm} is complete.

\end{document}